\documentclass[MSNbibl,number,citesort,seceqn,dvips]{arxbj}
\usepackage{dcolumn}
\usepackage{graphicx}

\aid{0}
\volume{20}
\issue{4}
\pubyear{2014}
\firstpage{2076}
\lastpage{2101}
\doi{10.3150/13-BEJ551} 

\makeatletter
\newcolumntype{d}[1]{D{.}{.}{#1}}
\newtheorem{theorem}{Theorem}[section]
\newtheorem{corollary}{Corollary}[section]
\newremark{remark}{Remark}[section]

\renewcommand{\citep}[1]{\cite{#1}}
\makeatother

\begin{document}
\begin{frontmatter}

\title{Stochastic monotonicity and continuity properties of functions
defined
on Crump--Mode--Jagers branching processes, with application to
vaccination in epidemic modelling}
\runtitle{Crump--Mode--Jagers branching processes}

\begin{aug}
\author[1]{\inits{F.}\fnms{Frank} \snm{Ball}\thanksref{1}\ead[label=e1]{frank.ball@nottingham.ac.uk}},
\author[2]{\inits{M.}\fnms{Miguel} \snm{Gonz\'alez}\thanksref{2,e2}\ead[label=e2,mark]{mvelasco@unex.es}},
\author[2]{\inits{R.}\fnms{Rodrigo}~\snm{Mart{\'i}nez}\corref{}\thanksref{2,e3}\ead[label=e3,mark]{rmartinez@unex.es}}
\and\break
\author[3]{\inits{M.}\fnms{Maroussia} \snm{Slavtchova-Bojkova}\thanksref{3}\ead[label=e4]{bojkova@fmi.uni-sofia.bg}}
\address[1]{School of Mathematical Sciences, The University of
Nottingham, Nottingham NG7 2RD, United Kingdom.
\printead{e1}}

\address[2]{Department of Mathematics, University of Extremadura, Avda,
Elvas s/n, 06071-Badajoz, Spain.
\printead{e2,e3}}

\address[3]{Faculty of Mathematics and Informatics, Sofia
University and Institute of Mathematics and Informatics, Bulgarian
Academy of Sciences, Bulgaria. \printead{e4}}
\end{aug}

\received{\smonth{6} \syear{2012}}
\revised{\smonth{7} \syear{2013}}

%
\begin{abstract}
This paper is concerned with Crump--Mode--Jagers branching
processes, describing spread of an epidemic depending on the
proportion of the population that is vaccinated. Births in the
branching process are aborted independently with a time-dependent
probability given by the fraction of the population vaccinated.
Stochastic monotonicity and continuity results for a wide class of
functions (e.g., extinction time and total number of births over all
time) defined on such
a branching process are proved using coupling arguments, leading
to optimal vaccination schemes to control corresponding functions
(e.g., duration and final size) of epidemic outbreaks. The theory is
illustrated by
applications to the control of the duration of mumps outbreaks in Bulgaria.
\end{abstract}

%
\begin{keyword}
\kwd{coupling}
\kwd{general branching process}
\kwd{Monte-Carlo method}
\kwd{mumps in Bulgaria}
\kwd{SIR epidemic model}
\kwd{time to extinction}
\kwd{vaccination policies}
\end{keyword}

\end{frontmatter}

\section{Introduction}

Branching processes have been applied widely to model epidemic
spread (see, e.g., the monographs by Andersson and Britton \citep
{hakanbritton},
Daley and Gani \citep{daleygani} and Mode and Sleeman \citep
{modesleemam}, and the review by
Pakes \citep{pakes}). The process describing the number of infectious
individuals in an epidemic model may be well approximated by a
branching process if the population is homogeneously mixing and
the number of infectious individuals is small in relation to the
total size of the susceptible population, since under these
circumstances the probability that an infectious contact is with a
previously infected individual is negligible (see, e.g.,
Isham \citep{isham}). Such an approximation dates back to the pioneering
works of Bartlett \citep{bartlett} and Kendall \citep{kendall}, and can
be made
mathematically precise by showing convergence of the epidemic
process to a limiting branching process as the number of
susceptibles tends to infinity (see Ball \citep{ball83}, Ball and
Donnelly \citep{ball}
and Metz \citep{metz}). The approximation may also be extended to
epidemics in populations that are not homogeneously mixing, for
example, those containing small mixing units such as households and
workplaces (see Pellis \textit{et al.} \citep{pellis}).

Before proceeding we give outline descriptions of some common branching
process models (see, e.g., Jagers \citep{jage} for further details), which
describe the evolution of a single-type population. In all of these models,
individuals have independent and identically distributed reproduction processes.
In a Bienaym{\'e}--Galton--Watson branching process, each individual
lives for
one unit of time and then has a random number of children, distributed
according to a random variable,
$\zeta$ say. In a Bellman--Harris branching process (BHBP), each
individual lives until a random age,
distributed according to a random variable $I$ say, and then has a
random number of children, distributed according to $\zeta$, where $I$
and $\zeta$ are independent. The
Sevast'yanov branching process (SBP) is defined similarly, except $I$
and $\zeta$ may be dependent, so the number of children an individual
has is correlated with that individual's lifetime. Finally,
in a general branching process, also called a Crump--Mode--Jagers (CMJ)
branching process, each individual lives until a random age,
distributed according to $I$, and reproduces at ages according to a
point process $\xi$. More precisely, if an individual, $i$ say having
reproduction variables $(I_i,\xi_i)$, is born at time $b_i$ and $0 \le
\tau_{i1} \le\tau_{i2}\le
\cdots\le I_i$ denote the points of $\xi_i$, then individual $i$ has
one child at each of times $b_i+\tau_{i1}, b_i+\tau_{i2}, \ldots\,$.

This paper is primarily concerned with models for epidemics of
diseases, such as measles, mumps and avian influenza, which follow
the so-called SIR (Susceptible $\to$ Infective $\to$ Removed)
scheme in a closed, homogeneously mixing population or some of its
extensions. A key epidemiological parameter for such an epidemic
model is the basic reproduction number $R_0$ (see
Heesterbeek and Dietz \citep{heesterbeek}), which in the present
setting is given by the
mean of the offspring distribution of the approximating branching
process. In particular a major outbreak (i.e., one whose size is
of the same order as the population size) occurs with nonzero
probability if and only if $R_0>1$. Suppose that $R_0>1$ and a
fraction $c$ of the population is vaccinated with a perfect
vaccine in advance of an epidemic. Then $R_0$ is reduced to
$(1-c)R_0$, since a proportion $c$ of infectious contacts is with
vaccinated individuals. It follows that a major outbreak is
almost surely prevented if and only if $c \ge1-R_0^{-1}$. This
well-known result, which gives the critical vaccination coverage
to prevent a major outbreak and goes back at least to 1964
(e.g., Smith \citep{Smith}), is widely used to inform public health authorities.

As a consequence of the above result, many analyses of vaccination
strategies in the epidemic modelling literature have focussed on
reducing $R_0$ to its critical value of one. However, if the
population is large, both the total size and the duration of an
outbreak may still be appreciable. Indeed, in the limit as the
population size tends to infinity, {when $R_0=1$,} both of these
quantities have infinite expectation under any plausible modelling
assumptions. In practice, there may be a cost associated with an
individual contracting the disease being modelled, in which case
it is of interest to determine vaccination strategies which reduce
the expected value of the total cost of an outbreak to an
acceptable level. Alternatively, it may be desired to control the
duration of an outbreak, for example, if the presence of an
outbreak means that restrictions are placed on the population
within which it is spreading. {Clearly, for large populations,
both of these aims necessitate that $R_0$ is reduced to somewhat
less than one.} The above remarks pertain to the common situation
of controlling an epidemic that is in its increasing phase. A
different situation arises with diseases, such as measles and
mumps, which are controlled by mass vaccination but small
outbreaks still occur among unvaccinated individuals.
Supplementary vaccination may be used to reduce the size or
duration of such outbreaks (as in the illustrative example of
mumps in Bulgaria in Section~\ref{illus} of this paper). A
similar phenomenon occurs with pathogens, such as monkeypox virus,
which primarily affect animals but spill over into human
populations giving stuttering chains of human-to-human
transmission (Lloyd-Smith \textit{et al.} \cite{lsmith}). In at least some
of the above scenarios, it may be the case that a specific
vaccination level cannot be achieved immediately but rather the
fraction of the population that is vaccinated will be
time-dependent. The aim of this paper is to develop a methodology
based on branching processes for addressing the above issues in a
unified fashion.

Gonz\'alez \textit{et al.} \citep{gmb1,gmb2} studied properties of the time
to extinction of an epidemic given that a fraction $c$ of
individuals is vaccinated, when the number of infectious
individuals in the population is modelled by a continuous-time
BHBP and a (more general) continuous-time SBP,
respectively. In an earlier work, De Serres \textit{et al.} \citep{dsgf} used a
discrete-time Bienaym{\'e}--Galton--Watson branching process to
study the spread of an infectious disease under various control
measures, specifically to estimate the effective
(i.e., post-control) value of $R_0$ from observations on size and
durations of small outbreaks.
The main objective in Gonz\'alez \textit{et al.} \citep{gmb1,gmb2} was to
determine the optimal proportion of susceptible individuals which
has to be vaccinated so that the mean (or given quantile of the)
extinction time of the disease is less than some specified value.
To that end, stochastic monotonicity and continuity properties of
the distribution function and mean of the time that the infection
survives, depending on the vaccination coverage rate were first
determined.

{In the present paper, we} extend the results in Gonz\'alez \textit{et
al.} \citep{gmb1,gmb2} in several directions that are both
practically and theoretically important. First, we assume that the
spread of infection is modelled as a CMJ branching process. The
CMJ branching process is appropriate for modelling the early
stages of a very wide variety of SIR epidemics, and includes both
BHBP and SBP as special cases. Second, we consider more general
vaccination processes. In Gonz\'alez \textit{et
al.} \citep{gmb1,gmb2} it was assumed that the fraction of
the population that is vaccinated remained constant with time. We
now allow this fraction to be an arbitrary but specified function
of time, thus capturing for example the setting in which people
are vaccinated as the disease spreads. Third, we consider the
control of more general functions of the epidemic process.
Gonz\'alez \textit{et al.} \citep{gmb1,gmb2} focused on
controlling the duration of the epidemic. The methods developed in
this paper are applicable to a wide class of functions of the
epidemic process. In addition to the duration of an outbreak, this
class includes, for example, the total number of people infected
and the maximum number of infected people present during the
epidemic.

The methodology of the paper is very different from that of
Gonz\'alez \textit{et al.} \citep{gmb1,gmb2}. The key stochastic
monotonicity and
continuity results in these papers were obtained by analysis
of integral equations governing properties of the time to
extinction of the branching process. In the present paper, a main
tool is coupling and, in particular, a pruning method of constructing a
realisation of a vaccinated process from that of the
corresponding unvaccinated process. As indicated in Section~\ref{conc},
this methodology is very powerful and applicable to a broad range
of processes.

The remainder of the paper is organised as follows. In
Section~\ref{modcoup}, we describe a very general model for an SIR
epidemic in a closed, homogeneously mixing community and explain
why its early spread may be approximated by a CMJ branching
process. We introduce a very general vaccination process and give
the basic coupling construction for obtaining a realisation of the
vaccinated epidemic process from that of the unvaccinated process.
The theoretical results of the paper are given in
Section~\ref{sectmoncty}. In Section~\ref{monprune}, we introduce
functions of a realisation of a CMJ branching process that are
monotonically decreasing with pruning. Examples of such functions
include the extinction time, the maximum population size over all
time and the total number of births over all time. Then we prove
in general, that is, independently of the function, monotonicity and
continuity properties of the mean (Section~\ref{monmean}),
distribution function (Section~\ref{mondsn}) and quantiles
(Section~\ref{monquantile}) of such functions. In
Section~\ref{optimal}, we use the previous results to define
optimal vaccination policies based on mean and quantiles. The
theory is then specialised in Section~\ref{timeextinction} to the
extinction time of an outbreak. The {methodology} is illustrated
in Section~\ref{illus} with applications to mumps in Bulgaria,
where vaccination is targeted at reducing the duration of an
outbreak. The paper ends with some concluding comments in
Section~\ref{conc}.

\section{Model and coupling construction}
\label{modcoup}

Consider first the following model for the spread of an epidemic
in a closed, homogeneously mixing population. Initially there are
$a$ infectives and $N$ susceptibles. Infectious individuals have
independent and identically distributed life histories
$\mathcal{H}=(I,\xi)$, where $I$ is the time elapsing between an
individual's infection and his/her eventual removal or death and
$\xi$ is a point process of times, relative to an individual's
infection, at which infectious contacts are made. Each contact is
with an individual chosen independently and uniformly from the
population. If a contact is with an individual who is susceptible,
then that individual becomes infected and itself makes contacts
according to its life history. If a contact is with an individual
who is not susceptible, then nothing happens. The epidemic ceases
as soon as there is no infective present in the population. Note
that, for simplicity, we assume that every infectious contact with
a susceptible necessarily leads to that susceptible becoming
infected. The model is easily extended to the situation when each
contact with a susceptible is successful (i.e., leads to infection)
independently with probability $p$ by letting
$\mathcal{H}=(I,\xi')$, where $\xi'$ is a suitable thinning of
$\xi$.

The above model is essentially that introduced by Ball and
Donnelly \citep{ball}, who noted that it included as special cases
a range of specific models that had hitherto received considerable
attention in the literature. For example, SIR and SEIR
(Susceptible $\to$ Exposed (i.e., latent) $\to$ Infective $\to$
Removed) models come under the above framework. The only
difference between the above model and that in Ball and
Donnelly \citep{ball} is that, in the latter, each contact is with
an individual chosen independently and uniformly from the $N$
initial susceptibles (rather than from the entire population of
$N+a$ individuals). In Ball and Donnelly \citep{ball}, a coupling
argument (which also holds for the present model) is used to prove
strong convergence, as the number of initial susceptibles $N \to
\infty$ (with the number of initial infectives $a$ held fixed), of
the process of infectives in the epidemic model to a CMJ branching
process (see Jagers \citep{jage}), in which a typical individual
lives until age $I$ and reproduces at ages according to $\xi$.
Thus for large $N$, the epidemic may be approximated by the CMJ
branching process. The approximation assumes that every contact
is with a susceptible individual. The proof in Ball and
Donnelly \citep{ball} may be extended to epidemics other than SIR,
for example, SIS (Susceptible $\to$ Infective $\to$ Susceptible) and SIRS
(Susceptible $\to$ Infective $\to$ Removed $\to$ Susceptible), by
suitably generalizing the life history $\mathcal{H}$ to allow for
removed individuals to become susceptible again (see,
e.g., Ball \citep{ball99} in the context of epidemics among a
population partitioned into households). Indeed, for a very broad
class of homogeneously mixing epidemic models, {that covers all of
the common stochastic formulations of infectious disease spread,}
the early stages of an epidemic in a large population with few
initial infectives may be approximated by a CMJ branching process.

This paper is concerned with the use of vaccination schemes to
control an epidemic, for example, in terms of its duration
or of the total number of individuals infected. We are
thus interested in the short-term behaviour of the epidemic, so we
model the epidemic as a CMJ branching process, $Z=\{Z(t)\dvtx t\geq
0\}$, where $Z(t)$ denotes the number of infected individuals at
time $t$. Thus $Z(0)$, which we assume to be fixed, represents the
number of infected individuals at the beginning of the outbreak.

We model the vaccination process by a function $\alpha\dvtx [0,\infty)
\to[0,1]$, such that $\alpha(t)$ is the proportion of the
population that are immune at time $t$ ($t\ge0$). Thus, the
probability that a contact at time $t$ is with a susceptible (i.e., non-immune)
individual is $1-\alpha(t)$. If the vaccine is perfect, that is, it
confers immunity immediately with probability one, then
$\alpha(t)$ is given by the proportion of the population that has
been vaccinated by time $t$. If the vaccine is imperfect then that
is implicitly included in the function $\alpha$. For example, if
the vaccine is all-or-nothing (i.e., it renders the vaccinee
completely immune with probability $\varepsilon$, otherwise it has
no effect), then $\alpha(t)=\varepsilon\tilde{\alpha}(t)$, where
$\tilde{\alpha}(t)$ is the proportion of the population that has
been vaccinated by time $t$. Note that if the immunity conferred
by vaccination does not wane then $\alpha$
is nondecreasing in $t$. We denote by
$Z_{\alpha}=\{Z_{\alpha}(t)\dvtx t\geq0\}$ the vaccination version of
$Z$, in which each birth in $Z$ is aborted independently, with
probability $\alpha(t)$ if the birth time is at time $t$.

Let $\mathcal{A}$ be the space of all functions $\alpha\dvtx [0,\infty)
\to[0,1]$. We construct coupled realizations of $Z$ and
$Z_{\alpha}$ $(\alpha\in\mathcal{A})$ on a common probability
space $(\Omega, \mathcal{F},P)$ as follows. Let $(\Omega_1,
\mathcal{F}_1,P_1)$ be a probability space on which are defined
independent life histories $\mathcal{H}_1,\mathcal{H}_2,\dots\,$,
each distributed as $\mathcal{H}$, which are pieced together in
the obvious fashion to construct a realization of $Z$. More
specifically, the life histories
$\mathcal{H}_1,\mathcal{H}_2,\dots,\mathcal{H}_a$ are assigned to
the $a$ initial infectives and, for $i=1,2,\dots\,$, the $i$th
individual born in $Z$ is assigned the life history
$\mathcal{H}_{a+i}$. Note that with this construction $Z$ may be
viewed as a tree, which is augmented with birth and death times of
branches. Let $(\Omega_2, \mathcal{F}_2,P_2)$ be a probability
space on which is defined a sequence $U_1,U_2,\dots$ of
independent random variables, each uniformly distributed on
$(0,1)$. Let $(\Omega, \mathcal{F},P)=(\Omega_1 \times\Omega_2,
\mathcal{F}_1 \times\mathcal{F}_2, P_1 \times P_2)$. Then, for
$\alpha\in\mathcal{A}$, a realization of $Z_{\alpha}$ is
constructed on $(\Omega, \mathcal{F},P)$ as follows. For
$i=1,2,\dots\,$, let $b_i$ denote the time of the $i$th birth in
$Z$, if such a birth occurs. Then this birth is deleted in
$Z_{\alpha}$ if and only if $U_i\leq\alpha(b_i)$. If a birth is
deleted in $Z_{\alpha}$, then none of the descendants of that
individual in $Z$ occurs in $Z_{\alpha}$. Thus, if the $j$th birth
in $Z$ is such a descendant then $U_j$ is redundant in the
construction of $Z_{\alpha}$. With the tree setting in mind, the
process of deleting an individual and all of its descendants is
called \emph{pruning}. For a previous use of pruning
in a branching process framework see, for example, Aldous and Pitman
\cite{aldous}.

Finally, we give some notation concerned with functions in $\mathcal
{A}$, which
will be used throughout the paper.
For $\alpha, \alpha' \in\mathcal{A}$, write $\alpha\prec
\alpha'$ if $\alpha(t) \le\alpha'(t)$ for all $t \in[0,\infty)$.
Also, for any $c \in[0,1]$ and any
$t_0 \ge0$, define the function $\alpha_c^{t_0} \in
\mathcal{A}$ by
\[
\alpha_c^{t_0}(t)= %
\cases{ 0 & \quad \mbox{if} $t<
t_0$,
\cr
c & \quad \mbox{if} $t\geq t_0$. } %
\]
Thus, for example, $\alpha_c^0$ denotes the constant function
equal to $c$ and $\alpha_0^0$ denotes the constant function equal
to $0$.

\section{Monotonicity and continuity properties depending on
vaccination function \texorpdfstring{$\alpha$}{alpha}} \label{sectmoncty}
\subsection{Functions \texorpdfstring{$f(Z_{\alpha})$}{f(Z alpha)} monotone to pruning}
\label{monprune}

Let $f(Z)$ be any nonnegative function of $Z$ taking values in
the extended real line $\mathbb{R} \cup\{\infty\}$ and, for
$\alpha\in\mathcal{A}$, let $\mu_{\alpha}^f=\mathrm{E}[f(Z_{\alpha
})]$. Again with the tree setting in mind, we say
that $f$ is monotonically decreasing with pruning, and write $f
\in\mathcal{P}$, if $f(Z^P) \le f(Z)$ almost surely
whenever $Z^P$ is obtained from $Z$ by pruning. For an event, $E$
say, let $1_E$ denote the indicator function of $E$. Examples of
functions that are monotonically decreasing with pruning include:
\begin{enumerate}[(iii)]
\item[(i)] the extinction time $T=\inf\{t\ge0\dvtx Z(t)=0\}$ and
$1_{\{T > t\}}$, where $t \in[0,\infty)$ is fixed;
\item[(ii)]
the maximum population size (number of infected individuals in
the epidemic context) over all time, $M=\sup_{t \ge0} Z(t)$
and $1_{\{M>x\}}$, where $x \in[0,\infty)$ is fixed;
\item
[(iii)] $N(t)$, the total number of births (new infections in the
epidemic context) in $(0,t]$, where $t \in[0,\infty)$ is fixed,
and the total number of births over all time (outbreak total size in
the epidemic
context) $N(\infty)=\lim_{t \to\infty}N(t)$, together with the
corresponding indicator functions $1_{\{N(t) > x\}}$ and
$1_{\{N(\infty) > x\}}$, where $x \in[0,\infty)$ is fixed.
\end{enumerate}

Throughout the paper, we assume that $Z$ is non-explosive,
that is, that $\mathrm{P}(N(t)<\infty)=1$ for any $t \in(0,\infty)$.
Conditions which guarantee this property may be found in
Jagers \citep{jage}, Section~6.2.

\subsection{Monotonicity and continuity of mean of
\texorpdfstring{$f(Z_{\alpha})$}{f(Z alpha)}} \label{monmean}

In this subsection, we derive monotonicity and continuity
properties of $\mathrm{E}[f(Z_{\alpha})]$, when viewed as a function
of the vaccination process $\alpha$, for functions $f$ that are monotonically
decreasing with pruning.
%
\begin{theorem}\label{teoa1} If $\alpha, \alpha' \in\mathcal{A}$
satisfy $\alpha\prec\alpha'$ and $f \in\mathcal{P}$, then $\mu
_{\alpha}^f \ge
\mu_{\alpha'}^f$.
\end{theorem}
\begin{pf}
The result follows immediately from the above construction of $Z$
and $Z_{\alpha}$, $\alpha\in\mathcal{A}$, on $(\Omega,
\mathcal{F},P)$, since $f$ is monotonically decreasing with
pruning and $Z_{\alpha'}$ may be obtained from $Z_{\alpha}$ by
successive prunings.
\end{pf}

We now give conditions under which $\mu_{\alpha}^f$ is continuous
in $\alpha$. For $\alpha,\alpha' \in\mathcal{A}$, let
$\|\alpha-\alpha'\|=\sup_{t \in
[0,\infty)}|\alpha(t)-\alpha'(t)|$ and, for $t>0$, let $\|\alpha-\alpha
'\|_t=\sup_{s \in
[0,t]}|\alpha(s)-\alpha'(s)|$. For $t>0$, write $f \in\mathcal{P}_t$
if $f \in
\mathcal{P}$ and $f(Z)$ depends on $Z$ only through $\{Z(s)\dvtx 0 \le
s \le t\}$. Let $m$ be the offspring mean for $Z$. For $c \in[0,1]$,
let $m_c$ denote the offspring mean of $Z_{\alpha_c^0}$, so
$m_c=(1-c)m$. Further, let $c_{\mathrm{inf}}=\max(0,1-m^{-1})$ and note that
$m_{c_{\mathrm{inf}}}\le1$. For $t_0 \ge0$ and $c \in[0,1]$, let
\begin{eqnarray*}
\mathcal{A}(c,t_0)=\bigl\{\alpha\in \mathcal{A}\dvtx  \alpha(t) \ge c
\mbox{ for all }t \ge t_0\bigr\}.
\end{eqnarray*}

\begin{theorem}\label{teoa2}
\begin{enumerate}[(b)]
\item[(a)] Fix $t>0$, let $f \in\mathcal{P}_t$ and
suppose that there exists a non-negative real-valued function $\hat{f}$,
with $\mathrm{E}[\hat{f}(Z)] < \infty$, such that, for
$P$-almost all $\omega\in\Omega$,
%
\begin{equation}
\label{Xbound2} f\bigl(Z_{\alpha}(\omega)\bigr) \le\hat{f}\bigl(Z(\omega)
\bigr)\qquad  \mbox{for all } \alpha\in \mathcal{A}.
\end{equation}
Then, for each $\varepsilon>0$, there exists
$\eta=\eta(\varepsilon)>0$ such that for all $\alpha, \alpha' \in
\mathcal{A}$ satisfying
$\|\alpha-\alpha'\|_t \le\eta$,
%
\begin{equation}
\label{mucontinuity} \bigl|\mu_{\alpha}^f-\mu_{\alpha'}^f\bigr|
\le\varepsilon.
\end{equation}
%
%
\item[(b)] Suppose that $m < \infty$. Let $f \in\mathcal{P}$ and
$t_0\geq
0$, and suppose that there exists a non-negative real-valued
function $\hat{f}(Z_{\alpha_{c_{\mathrm{inf}}}^{t_0}})$,
with $\mathrm{E}[\hat{f}(Z_{\alpha_{c_{\mathrm{inf}}}^{t_0}})] < \infty
$, such that, for
$P$-almost all $\omega\in\Omega$,
%
\begin{equation}
\label{Xbound1} f\bigl(Z_{\alpha}(\omega)\bigr) \le \hat{f}
\bigl(Z_{\alpha_{c_{\mathrm{inf}}}^{t_0}}(\omega)\bigr)\qquad  \mbox{for all } \alpha\in
\mathcal{A}(c_{\mathrm{inf}}, t_0).
\end{equation}
Then, for each $\varepsilon>0$, there exists
$\eta=\eta(\varepsilon)>0$ such that (\ref{mucontinuity}) holds
for all $\alpha, \alpha' \in
\mathcal{A}(c_{\mathrm{inf}}, t_0)$ satisfying
$\|\alpha-\alpha'\| \le\eta$.
\end{enumerate}
\end{theorem}
\begin{pf}
(a)
For $n=1,2,\dots$ and $\alpha, \alpha' \in\mathcal{A}$, let
\begin{eqnarray*}
B_n\bigl(\alpha,\alpha'\bigr)= \bigcap
_{i=1}^n \bigl\{\omega\in \Omega\dvtx U_i(
\omega) \notin \bigl(\min\bigl(\alpha(b_i),\alpha'(b_i)
\bigr),\max\bigl(\alpha(b_i),\alpha'(b_i)
\bigr)\bigr]\bigr\},
\end{eqnarray*}
and let $B_0(\alpha,\alpha')=\Omega$. Now $\mathrm{P}(N(t)<\infty)=1$,
since $Z$ is non-explosive. Observe that if $\omega
\in B_{N(t)}(\alpha,\alpha')$ then, by construction,
$Z_{\alpha}(s,\omega)=Z_{\alpha'}(s,\omega)$ for all $s \in[0,t]$,
whence $f(Z_{\alpha}(\omega))=f(Z_{\alpha'}(\omega))$ since $f \in
\mathcal{P}_t$. Now, for any $\alpha
\in\mathcal{A}$,
\begin{eqnarray*}
\mu_{\alpha}^f=\mathrm{E} \bigl[f(Z_{\alpha})1_{B_{N(t)}(\alpha,\alpha
')}
\bigr]+\mathrm{E} \bigl[f(Z_{\alpha}) 1_{B_{N(t)}^c(\alpha,\alpha')} \bigr],
\end{eqnarray*}
where $B_{N(t)}^c(\alpha,\alpha')=\Omega\setminus B_{N(t)}(\alpha
,\alpha')$.
Thus, for any $\alpha, \alpha' \in\mathcal{A}$,
\begin{eqnarray*}
\mu_{\alpha}^f-\mu_{\alpha'}^f=\mathrm{E}
\bigl[f(Z_{\alpha}) 1_{B_{N(t)}^c(\alpha,\alpha')} \bigr]-\mathrm{E} \bigl[f(Z_{\alpha'})
1_{B_{N(t)}^c(\alpha,\alpha')} \bigr],
\end{eqnarray*}
whence, since $f$ is nonnegative,
\begin{eqnarray*}
\bigl|\mu_{\alpha}^f-\mu_{\alpha'}^f\bigr| \le
\mathrm{E} \bigl[\hat{f}(Z) 1_{B_{N(t)}^c(\alpha,\alpha')} \bigr].
\end{eqnarray*}
Now
\begin{eqnarray*}
\mathrm{E} \bigl[\hat{f}(Z) 1_{B_{N(t)}^c(\alpha,\alpha')} \bigr]= \mathrm{E} \bigl[\hat{f}(Z)
\mathrm{E} [1_{B_{N(t)}^c(\alpha,\alpha
')}|Z ] \bigr].
\end{eqnarray*}
Further, (i) $Z$ determines $N(t)$ and (ii) $(U_1,U_2,\dots)$ is
independent of $Z$,
so, $P$-almost surely,
\begin{eqnarray*}
\mathrm{E} [1_{B_{N(t)}^c(\alpha,\alpha')}|Z ]&=& 1-\prod_{i=1}^{N(t)}
\bigl(1-\bigl|\alpha(b_i)-\alpha'(b_i)\bigr|\bigr)
\\
&\le&1-(1-\delta)^{N(t)},
\end{eqnarray*}
where $\delta=\|\alpha-\alpha'\|_t$. Hence, $P$-almost surely,
\begin{eqnarray*}
\mathrm{E} [1_{B_{N(t)}^c(\alpha,\alpha')}|Z ] \le \mathrm{E} [1_{B_{N(t)}^c(\alpha_0^0,\alpha_{\delta}^0)}|Z ],
\end{eqnarray*}
whence, for $\alpha, \alpha' \in\mathcal{A}$,
%
\begin{eqnarray}
\label{hatmubound} \bigl|\mu_{\alpha}^f-\mu_{\alpha'}^f\bigr|
&\le& \mathrm{E} \bigl[\hat {f}(Z)1_{B_{N(t)}^c(\alpha_0^0,\alpha_{\delta}^0)} \bigr]
\nonumber
\\[-8pt]\\[-8pt]
&=&\hat{\mu}_t(\delta)\qquad  \mbox{say}.\nonumber
\end{eqnarray}

Now $\mathrm{P}(N(t)<\infty)=1$, so $P$-almost surely,
\begin{eqnarray*}
\hat{f}(Z)1_{B_{N(t)}^c(\alpha_0^0,\alpha_{\delta}^0)} \to0 \qquad \mbox{as } \delta\downarrow0
\end{eqnarray*}
(in fact $\hat{f}(Z)1_{B_{N(t)}^c(\alpha_0^0,\alpha_{\delta}^0)}=0$ for
all $\delta\in[0,\delta^*)$, where $\delta^*= \min(U_1,U_2,\dots,
U_{N(t)})$), so by the dominated convergence theorem $\hat{\mu}_t(\delta
) \to0$ as $\delta\downarrow0$. Thus, given $\varepsilon>
0$, there exists $\eta$ such that $\hat{\mu}_t(\delta) \le\varepsilon$
for all $\delta\in(0,\eta)$ and the theorem follows using (\ref{hatmubound}).


(b) For $\alpha\in\mathcal{A}(c_{\mathrm{inf}}, t_0)$, the
process $Z_\alpha$ can be viewed as a vaccinated version of the process
$Z_{\alpha_{c_{\mathrm{inf}}}^{t_0}}$ with vaccination function $\tilde
{\alpha}$ given by
\[
\tilde{\alpha}(t)= %
\cases{ \alpha(t) &\quad  \mbox{if} $t<
t_0$,
\cr
\displaystyle \frac{\alpha
(t)}{1-c_{\mathrm{inf}}} & \quad \mbox{if} $t\geq t_0$. }
\]
Note that $Z_{\alpha_{c_{\mathrm{inf}}}^{t_0}}$ has
offspring mean $m$ until time $t_0$, and $m_{c_{\mathrm{inf}}}\le1$
after time $t_0$. Thus, since $Z$ is non-explosive (so $\mathrm
{P}(Z(t_0)<\infty)=1$), the total
number of births over all time in $Z_{\alpha_{c_{\mathrm{inf}}}^{t_0}}$
(i.e., $N_{\alpha_{c_{\mathrm{inf}}}^{t_0}}(\infty)$) is finite almost
surely. Also, $\|\tilde{\alpha}-\tilde{\alpha}'\| \le
(1-c_{\mathrm{inf}})^{-1} \|\alpha-\alpha'\|$. The proof then proceeds
as in part (a), but with
$Z$ and $N(t)$ replaced by $Z_{\alpha_{c_{\mathrm{inf}}}^{t_0}}$ and
$N_{\alpha_{c_{\mathrm{inf}}}^{t_0}}(\infty)$, respectively, and $\alpha
, \alpha'$ replaced by $\tilde{\alpha}, \tilde{\alpha}'$.
\end{pf}
%
%
\begin{remark}\label{rr1}
\begin{enumerate}[(b)]
\item[(a)] Suppose that $m \le1$. Then $c_{\mathrm{inf}}=0$ and it
follows that $Z_{\alpha_{c_{\mathrm{inf}}}^{t_0}}=Z$ and $\mathcal
{A}(c_{\mathrm{inf}}, t_0)=\mathcal{A}$. Thus, for any $f \in\mathcal
{P}$, Theorem~\ref{teoa2}(b) implies that, for any $\varepsilon>0$,
there exists
$\eta=\eta(\varepsilon)>0$ such that (\ref{mucontinuity}) holds
for all $\alpha, \alpha' \in
\mathcal{A}$ satisfying
$\|\alpha-\alpha'\| \le\eta$.

\item[(b)] Suppose that $m > 1$ and $f \in\mathcal{P}$. Then the
argument used to prove Theorem~\ref{teoa2}(b) breaks down since $\mathrm
{P}(Z(\infty)<\infty)<1$. Thus, with our argument we can prove
continuity in $\alpha$ of $\mu_\alpha^f$ for $f \in\mathcal{P}_t$, for
any $t>0$, but not for $f \in\mathcal{P}$. However, this is no
restriction from a practical viewpoint since $t$ in Theorem~\ref{teoa2}(a),
or $t_0$ in Theorem~\ref{teoa2}(b), can be made arbitrarily large. For
example, in any real life-setting there will be a maximum time frame
over which it is of interest to evaluate the performance of a
vaccination process and $t$ or $t_0$ can be chosen accordingly.
\end{enumerate}
\end{remark}

\subsection{Monotonicity and continuity of distribution function
of~\texorpdfstring{$f(Z_{\alpha})$}{f(Z alpha)}} \label{mondsn}

Using the previous results, we
establish in this subsection monotonicity and continuity
properties of the distribution function of $f(Z_{\alpha})$. For $f
\in\mathcal{P}$ and $\alpha\in\mathcal{A}$, let
\begin{eqnarray*}
v_{\alpha}^f(x)=\mathrm{P}\bigl(f(Z_{\alpha}) \le x
\bigr)=1-\mathrm{E}[1_{\{f(Z_{\alpha})>x\}}], \qquad x \ge0,
\end{eqnarray*}
be the distribution function of the random variable
$f(Z_{\alpha})$.

For $\alpha\in\mathcal{A}$ and $t \in[0,\infty]$, let
$\phi_{N_\alpha(t)}(s)=\mathrm{E}[s^{N_\alpha(t)}]$ $(0 \le s \le1)$
denote the probability generating function of $N_\alpha(t)$.
Suppose that $P(N_\alpha(t)<\infty)=1$. Then
$\phi_{N_\alpha(t)}(1-)=1$ and $\phi_{N_\alpha(t)}^{-1}(u)$ is
well defined for all $u \in[u_{\alpha,t},1]$, where
$u_{\alpha,t}=\mathrm{P}(N_\alpha(t)=0)$. Extend the domain of
$\phi_{N_\alpha(t)}^{-1}$ by defining
$\phi_{N_\alpha(t)}^{-1}(u)=0$ for $u \in[0,u_{\alpha,t})$.
Define the function $\delta_{\alpha,t}\dvtx [0,1] \to[0,1]$ by
%
\begin{equation}
\label{deltaepsilon} \delta_{\alpha,t}(\varepsilon)=1-\phi_{N_\alpha(t)}^{-1}(1-
\varepsilon),\qquad  0 \le\varepsilon\le1.
\end{equation}
Note that $\delta_{\alpha,t}(\varepsilon)>0$ if $\varepsilon>0$
and $\lim_{\varepsilon\downarrow0}
\delta_{\alpha,t}(\varepsilon)=0$.
%
\begin{theorem}\label{teoa4}
\begin{enumerate}[(b)]
\item[(a)] Suppose that $f \in\mathcal{P}$ and $\alpha, \alpha'
\in\mathcal{A}$ satisfy $\alpha\prec\alpha'$.
Then
%
\begin{equation}
\label{valphainequ} v_{\alpha}^f(x) \le v_{\alpha'}^f(x)
\qquad \mbox{for all } 0 \le x \le\infty.
\end{equation}
\item[(b)] Fix $t>0$ and suppose that $f \in\mathcal{P}_t$. Then, for
any $\varepsilon>0$,
%
\begin{equation}
\label{nucontinuity} \sup_{0\leq x < \infty}\bigl|v_\alpha^f(x)-v_{\alpha'}^f(x)\bigr|
\leq\varepsilon
\end{equation}
for all $\alpha, \alpha' \in\mathcal{A}$ satisfying
$\|\alpha-\alpha'\|_t \le\delta_{\alpha_0^0,t}(\varepsilon)$.
\item[(c)] Suppose that $f \in\mathcal{P}$. Then, for any
$\varepsilon>0$, (\ref{nucontinuity}) holds for all $\alpha,
\alpha' \in\mathcal{A}(c_{\mathrm{inf}},t_0)$ satisfying
$\|\alpha-\alpha'\| \le
\delta_{\alpha_{c_{\mathrm{inf}}}^{t_0},\infty}(\varepsilon)$.
\end{enumerate}
\end{theorem}
\begin{pf}
(a) Fix $x \in[0,\infty)$ and let $\tilde{f}_x$ be the
function of $Z$ given by $\tilde{f}_x(Z)=1_{\{f(Z)>x\}}$. Then
$\tilde{f}_x \in\mathcal{P}$ and (\ref{valphainequ}) follows from
Theorem~\ref{teoa1}, since $v_{\alpha}^f(x)=1-\mathrm{E}[\tilde
{f}_x(Z_{\alpha})]$.

(b) For each $x \in[0,\infty)$,
\begin{eqnarray*}
\bigl|v_\alpha^f(x)-v_{\alpha'}^f(x)\bigr|=\bigl|
\mathrm{E}\bigl[\tilde{f}_x(Z_{\alpha
})\bigr]-\mathrm{E}\bigl[
\tilde{f}_x(Z_{\alpha'})\bigr]\bigr|
\end{eqnarray*}
and $\tilde{f}_x(Z_{\alpha}(\omega)) \le1$ for all $\alpha\in
\mathcal{A}$ and all $\omega\in\Omega$. Fix $t > 0$ and note that
$\tilde{f}_x \in\mathcal{P}_t$, since $f \in\mathcal{P}_t$ . It then
follows from
(\ref{hatmubound}), taking $\hat{f}(Z)=1$, that,
for $x \in[0,\infty)$ and $\alpha,\alpha' \in\mathcal{A}$,
%
\begin{equation}
\label{modvalpha} \bigl|v_\alpha^f(x)-v_{\alpha'}^f(x)\bigr|
\le \hat{\mu}_t\bigl(\bigl\|\alpha-\alpha'\bigr\|_t
\bigr),
\end{equation}
where, for $\delta\in[0,1]$,
\begin{eqnarray*}
\hat{\mu}_t(\delta)=\mathrm{P} \bigl(B_{N(t)}^c
\bigl(\alpha_0^0,\alpha_{\delta
}^0\bigr)
\bigr)=1-\mathrm{E} \bigl[(1-\delta)^{N(t)} \bigr]=1-\phi
_{N(t)}(1-\delta).
\end{eqnarray*}
Recall that $N(t)=N_{\alpha_0^0}(t)$ and note that
$P(N_{\alpha_0^0}(t)<\infty)=1$ since $Z$ is non-explosive. Thus,
$\phi_{N_{\alpha_0^0}(t)}^{-1}(u)$ is well defined for all $u \in
[0,1]$ and,  since
$1-\phi_{N_{\alpha_0^0}(t)}(1-\delta_{\alpha_0^0,t}(\varepsilon))\le
\varepsilon$, the theorem follows.

(c) The proof is similar to part (b) but with
$N_{\alpha_0^0}(t)$ replaced by $N_{\alpha_{c_{\mathrm
{inf}}}^{t_0}}(\infty)$.
\end{pf}
%
\begin{remark} \label{rr2}
\begin{enumerate}[(b)]
\item[(a)] Observe that the function $\delta_{\alpha_0^0,t}$,
defined using (\ref{deltaepsilon}), is independent of both $f$ and
$x$, so the uniform continuity of $v_{\alpha}^f(x)$, with respect
to $\alpha$, holds uniformly over all $f \in\mathcal{P}$ and all
$x \in[0,\infty)$.

\item[(b)] Similar to Remark~\ref{rr1}(a),
Theorem~\ref{teoa4}(c) shows that if $m \le1$ (so
$P(N(\infty)<\infty)=1$) and $f \in\mathcal{P}$ then, for any
$\varepsilon>0$, (\ref{nucontinuity}) holds for all $\alpha,
\alpha' \in\mathcal{A}$ satisfying $\|\alpha-\alpha'\| \le
\delta_{\alpha_0^0,\infty}(\varepsilon)$.
\end{enumerate}
\end{remark}

\subsection{Monotonicity and continuity of quantiles of
\texorpdfstring{$f(Z_{\alpha})$}{f(Z alpha)}} \label{monquantile}

In applications, we wish to
control the quantiles of $f(Z_{\alpha})$, so we now derive related
monotonicity and continuity properties. Fix $f \in\mathcal{P}$
and $\alpha\in\mathcal{A}$, and define, for $0<p<1$,
\begin{eqnarray*}
x_{\alpha,p}^f=\inf\bigl\{x\dvtx v_{\alpha}^f(x)\ge
p\bigr\},
\end{eqnarray*}
with the convention that $x_{\alpha,p}^f=\infty$ if
$v_{\alpha}^f(x)< p$ for all $x \in[0,\infty)$.
Thus,
$x_{\alpha,p}^f$ is the quantile of order $p$ of the random
variable $f(Z_{\alpha})$. For $\alpha\in\mathcal{A}$, let
$\mathcal{A}^+(\alpha)=\{\alpha' \in\mathcal{A}\dvtx \alpha\prec\alpha'\}$.
For a sequence
$\{\alpha_n\}$ and
$\alpha$ in $\mathcal{A}$, we define $\lim_{n \to\infty} \alpha
_n=\alpha$ to mean
$\lim_{n \to\infty}\|\alpha_n-\alpha\|=0$.
%
\begin{theorem}\label{teoa5}
Suppose that $f \in\mathcal{P}$ and $p \in(0,1)$.
\begin{enumerate}[(b)]
\item[(a)] If $\alpha,\alpha' \in\mathcal{A}$ satisfy $\alpha
\prec\alpha'$, then $x_{\alpha',p}^f \le x_{\alpha,p}^f$.
\item
[(b)] Suppose further that $f \in\mathcal{P}_t$ for some $t>0$ and
$\alpha\in\mathcal{A}$ is such that $x_{\alpha,p}^f < \infty$. Let
$\{\alpha_n\}$ be any sequence in $\mathcal{A}$
satisfying $\lim_{n \to\infty}
\alpha_n=\alpha$. Then
$\lim_{n \to\infty}
x_{\alpha_n,p}^f=x_{\alpha,p}^f$ in each of
the following cases:
\begin{enumerate}[(ii)]
\item[(i)] $\alpha_n \in
\mathcal{A}^+(\alpha)$ for all $n$;
\item[(ii)] $v_{\alpha}^f$ is continuous and
strictly increasing at $x_{\alpha,p}^f$.
\end{enumerate}
\end{enumerate}
\end{theorem}
\begin{pf}
(a) By Theorem~\ref{teoa4}(a), $\{x\dvtx v_{\alpha}^f(x)\ge p\}
\subseteq\{x\dvtx v_{\alpha'}^f(x)\ge p\}$, which implies
$x_{\alpha',p}^f \le x_{\alpha,p}^f$.

(b) Choose $t>0$ such that $f \in\mathcal{P}_t$. Suppose that
(i) holds. Let
$x_{\sup}=\limsup_{n \to\infty} x_{\alpha_n,p}^f$ and
$x_{\mathrm{inf}}=\liminf_{n \to\infty} x_{\alpha_n,p}^f$. Then by part (a),
$x_{\sup} \le x_{\alpha,p}^f$. Fix
$\varepsilon>0$. Then, since $\lim_{n \to\infty} \alpha_n=\alpha$
and $\|\alpha_n-\alpha\|_t\le\|\alpha_n-\alpha\|$, there exists
$n_0$ such that
$\|\alpha_n-\alpha\|_t\le\delta_{\alpha_0^0,t}(\varepsilon)$ for
all $n \ge n_0$, where $\delta_{\alpha_0^0,t}(\varepsilon)$ is
defined at (\ref{deltaepsilon}) -- recall that
$N(t)=N_{\alpha_0^0}(t)$. Now, $\alpha\prec\alpha_n$, hence, by
Theorem~\ref{teoa4}(a) and (b),
$v_{\alpha_n}^f(x)-v_{\alpha}^f(x)\le\varepsilon$, for all $x\ge
0$ and for all $n \ge n_0$. In particular, setting
$x=x_{\alpha_n,p}^f$ and noting that
$v_{\alpha_n}^f(x_{\alpha_n,p}^f) \ge p$ since $v_{\alpha_n}^f$ is
right-continuous, yields that $v_{\alpha}^f(x_{\alpha_n,p}^f) \ge
p-\varepsilon$ for all $n \ge n_0$. Hence, $v_{\alpha}^f(x_{\mathrm{inf}})
\ge p-\varepsilon$, since $v_{\alpha}^f$ is increasing and
right-continuous. This holds for all $\varepsilon>0$, so
$v_{\alpha}^f(x_{\mathrm{inf}}) \ge p$, whence $x_{\mathrm{inf}}\ge
x_{\alpha,p}^f$. Thus, $x_{\mathrm{inf}}=x_{\sup}=x_{\alpha,p}^f$, so
$\lim_{n \to\infty} x_{\alpha_n,p}^f=x_{\alpha,p}^f$, as
required.

Suppose that (ii) holds. First, we assume that $\alpha_n \prec
\alpha$ for all $n$. Then, by part (a), $x_{\mathrm{inf}} \ge
x_{\alpha,p}^f$. Note that $v_{\alpha}^f(x_{\alpha,p}^f)=p$, since
$v_{\alpha}^f$ is continuous at $x_{\alpha,p}^f$, and
$v_{\alpha}^f(x)>p$ for all $x>x_{\alpha,p}^f$, since
$v_{\alpha}^f$ is strictly increasing at $x_{\alpha,p}^f$. Fix
$x>x_{\alpha,p}^f$ and let $\varepsilon=v_{\alpha}^f(x)-p$, so
$\varepsilon>0$. As before, there exists $n_0$ such that
$\|\alpha_n-\alpha\|_t\le\delta_{\alpha_0^0,t}(\varepsilon)$ for
all $n \ge n_0$. It then follows from Theorem~\ref{teoa4} that
\begin{eqnarray*}
v_{\alpha}^f(x)-v_{\alpha_n}^f(x) \le
\varepsilon=v_{\alpha}^f(x)-p \qquad \mbox{for all }n \ge
n_0.
\end{eqnarray*}
Thus $v_{\alpha_n}^f(x) \ge p$ for all
$n \ge n_0$, whence
$x_{\alpha_n,p}^f \le x$ for all
$n \ge n_0$, which implies
that $x_{\sup} \le x$. Since this holds for any
$x>x_{\alpha,p}^f$, it follows that $x_{\sup} \le x_{\alpha,p}^f$,
which combined with $x_{\mathrm{inf}} \ge x_{\alpha,p}^f$ yields the
required result.

Now, we consider an arbitrary sequence
$\{\alpha_n\}$ that converges to $\alpha$. For $q=1,2,\dots\,$, define
functions $\alpha^+_q$
and $\alpha^-_q$ by
$\alpha^+_q(s)=\min\{\alpha(s)+\frac{1}{q},1\}$ and
$\alpha^-_q(s)=\max\{\alpha(s)-\frac{1}{q},0\}$ $(s \ge0)$.
Then $\lim_{q \to\infty} \alpha^+_q=\lim_{q \to\infty} \alpha
^-_q=\alpha$.
Further, $\alpha^-_q \prec\alpha\prec\alpha^+_q$ for each
$q=1,2,\dots\,$.
Hence, by part (i) and the above, $\lim_{q \to\infty} x_{\alpha
^+_q,p}^f=\lim_{q \to
\infty} x_{\alpha^-_q,p}^f=x_{\alpha,p}^f$. For any fixed $q \in
\mathbb{N}$, $\alpha_n \prec\alpha^+_q$ for
all sufficiently large $n$, so Theorem~\ref{teoa5}(a) implies that
$\liminf_{n \to\infty}
x_{\alpha_n,p}^f \ge x_{\alpha^+_q,p}^f$.
Letting $q \to\infty$ then yields that $ x_{\mathrm{inf}} \ge
x_{\alpha,p}^f$. A similar argument using the sequence
$\{\alpha^-_q\}$ shows that $x_{\sup} \le x_{\alpha,p}^f$, whence
$\lim_{n \to\infty}
x_{\alpha_n,p}^f= x_{\alpha,p}^f$, as
required.
\end{pf}
%
\begin{remark}\label{rr3}
\begin{enumerate}[(b)]
\item[(a)] It is straightforward to extend Theorem~\ref{teoa5}(b) to
a family of vaccination processes with a continuous index set, for
example, $\{\alpha_s:s \in\mathcal{I}\}$, where $\mathcal{I}$ is a
connected subset of $\mathbb{R}^d$ for some $d \in\mathbb{N}$.
Theorem~\ref{teoa5}(b) implies that, under appropriate conditions,
$\lim_{s \to s^*}{x_{\alpha_s,p}^f}= x_{\alpha_{s^*},p}^f$. We
use this extension when studying optimal vaccination policies in
the next subsection.

\item[(b)] Invoking Remark~\ref{rr2}(b) shows that if $m \le1$ then
Theorem~\ref{teoa5}(b) holds with $\mathcal{P}_t$ replaced by $\mathcal{P}$.
\end{enumerate}
\end{remark}

\subsection{Optimal vaccination policies based on mean and
quantiles} \label{optimal}

From the above monotonicity and
continuity properties of mean and quantiles, we propose next how
to choose optimal $\alpha$s, that is, optimal vaccination policies in
a sense that is made clear below, from a subset $\mathcal{A}^*$ of
$\mathcal{A}$. Fix $f\in\mathcal{P}$, $b> 0$ and $0<p<1$, and let
$\mathcal{A}_b^f=\{\alpha\in\mathcal{A}^*\dvtx  \mu_{\alpha}^f\leq
b\}$ and $\mathcal{A}_{p,b}^f=\{\alpha\in\mathcal{A}^*\dvtx  x_{\alpha
,p}^f\leq b\}$. {Notice that if, for example, $f$ is the
time to extinction, then $\mathcal{A}_b^f$ and
$\mathcal{A}_{p,b}^f$ comprise those vaccination policies in
$\mathcal{A}^*$ for which the mean and the quantile of order $p$,
respectively, of the time to extinction is less than or equal to
some bound $b$. Then it is of interest to search for optimal
vaccination policies which satisfy these properties.}

Then, if they exist, optimal vaccination policies based on the
mean are
\[
\mathop{\operatorname{argmax}}_{\alpha\in\mathcal{A}_b^f} \mu _{\alpha}^f
\]
and optimal vaccination policies based on the quantiles are
\[
\mathop{\operatorname{argmax}}_{\alpha\in\mathcal{A}_{p,b}^f} x_{\alpha,p}^f.
\]
We notice that the sets $\mathcal{A}_b^f$ and
$\mathcal{A}_{p,b}^f$ can be empty. If they are not empty, optimal
vaccination policies may not be unique when a total order is not
defined on the sets $\mathcal{A}_b^f$ and $\mathcal{A}_{p,b}^f$.
Otherwise, provided the conditions of
Theorems \ref{teoa1}, \ref{teoa2} and \ref{teoa5} are {satisfied},
the monotonicity and continuity properties of mean and quantiles
of $f(Z_{\alpha})$ proved in those theorems imply that there exist
unique $\alpha_{\mathrm{opt},b}^f\in\mathcal{A}_b^f$ and
$\alpha_{\mathrm{opt},p,b}^f\in\mathcal{A}_{p,b}^f$ such that
\[
\mu_{\alpha_{\mathrm{opt},b}^f}^f=\max_{\alpha\in\mathcal{A}_b^f}
\mu_{\alpha}^f \quad \mbox{and}\quad  x_{\alpha_{\mathrm{opt},p,b}^f,p}^f=\max
_{\alpha\in\mathcal{A}_{p,b}^f} x_{\alpha,p}^f.
\]
Intuitively, $\alpha_{\mathrm{opt},b}^f$ and
$\alpha_{\mathrm{opt},p,b}^f$ are the smallest vaccination policies in
$\mathcal{A}^*$ such that the mean\vspace*{2pt} and the $p$th quantile,
respectively, of $f(Z_{\alpha_{\mathrm{opt},b}^f})$ and
$f(Z_{\alpha_{\mathrm{opt},p,b}^f})$ are less than or equal to $b$. Before
giving some simple examples of $\mathcal{A}^*$, we discuss briefly
conditions that ensure the existence and uniqueness of optimal
policies.

For fixed $f \in\mathcal{P}$, define the binary relation
$\prec_f$ on $\mathcal{A}$ by $\alpha\prec_f \alpha'$ if and only
if $\mu_\alpha^f \le\mu_{\alpha'}^f$. Observe that, if $\alpha
\prec\alpha'$ then, by Theorem~\ref{teoa1}, {$\alpha' \prec_f
\alpha$} for any $f \in\mathcal{P}$. The relation $\prec_f$ is
not an ordering, because {$\alpha\prec_f \alpha'$ and $\alpha'
\prec_f \alpha$} imply only that $\mu_\alpha^f = \mu_{\alpha'}^f$
(and not that $\alpha=\alpha'$). However, we can consider the
equivalence relation $\sim_f$ on $\mathcal{A}$ defined by $\alpha
\sim_f \alpha'$ if and only if $\mu_\alpha^f = \mu_{\alpha'}^f$.
Then $\prec_f$ is a total ordering on the quotient set
$\mathcal{A}/\sim_f$, that is, the set of all possible equivalence
classes, using the obvious definition of $\prec_f$ on
$\mathcal{A}/\sim_f$.

Given a subset $\mathcal{A}^*$ of $\mathcal{A}$, a simple
condition that ensures the existence of ${\operatorname{argmax }}_{\alpha\in
\mathcal{A}_b^f}\mu_{\alpha}^f$ for any
fixed $b>0$ is that the set of real numbers
$\{\mu_{\alpha}^f\dvtx \alpha\in\mathcal{A}^*\}$ is closed. More
precisely, this ensures the existence of an equivalence class on
which the maximum is attained. To obtain a unique maximum
requires that $\prec_f$ is a total ordering on $\mathcal{A}^*$ (or
at least on $\mathcal{A}_b^f$ for fixed $b$). Note that even if
$\prec$ is a total ordering {on} $\mathcal{A}^*$,
Theorem~\ref{teoa1} does not ensure that $\prec_f$ is a total
ordering on $\mathcal{A}^*$. For the latter, we require that
$\mu_\alpha^f>\mu_{\alpha'}^f$ for all $\alpha, \alpha' \in
\mathcal{A}^*$ satisfying $\alpha\prec\alpha'$ and $\alpha\ne
\alpha'$. The coupling argument in Section~\ref{modcoup} can be
used to show that this holds for any practically useful $f$ and it
is assumed implicitly in the sequel. Similar arguments to the
above pertain for optimal vaccination policies based on quantiles.

A simple example of $\mathcal{A}^*$ is the set of constant
functions, that is, $\mathcal{A}^*=\{\alpha_c^0 \dvtx  0\leq c \leq1\}$.
On this set, the total order is defined by the order of the real
numbers. Another example is the set
$\mathcal{A}^*=\{\alpha_{M,t_v,p_0}\dvtx  M\geq0, 0\leq p_0\leq1,
0\leq t_v\leq{p_0^{-1}}\}$, where, for $s\geq0$,
%
\begin{equation}
\label{alphapara} \alpha_{M,t_v,p_0}(s)= %
\cases{ 0 & \quad \mbox{if }$s
\leq M$,
\cr
p_0(s-M) & \quad \mbox{if} $M <s\leq M+t_v$,
\cr
t_vp_0 & \quad \mbox{if} $M+t_v<s$. }
\end{equation}

For fixed $M$, $t_v$ and $p_0$, the function $\alpha_{M,t_v,p_0}$
describes the proportion of immune individuals in the population
when the vaccination process starts at time $M$, takes $t_v$ time
units and the proportion of individuals vaccinated per unit time
is $p_0$. We notice that a total order on $\mathcal{A}^*$ is not
possible. However, in practice, $M$ and $p_0$ are usually known
before vaccination begins, and therefore, the functions can be
parameterized through $t_v$ alone. For fixed $M$ and $p_0$,
denote $\alpha_{t_v}= \alpha_{M,t_v,p_0}$ and
$\mathcal{A}^*=\{\alpha_{t_v}\dvtx  {c_{\mathrm{inf}}p_0^{-1}} \le t_v \le
p_0^{-1}\}$. Then $\prec_f$ is a total ordering on $\mathcal{A}^*$
and Theorem~\ref{teoa2}(b) ensures that $\{\mu_{\alpha}^f\dvtx \alpha
\in\mathcal{A}^*\}$ is closed, so, provided $\mathcal{A}_b^f$ is
non-empty, the optimal vaccination policy exists and is unique.
Moreover, it and the corresponding optimal policies based on {the
mean and} quantiles are given by $\alpha_{t_{\mathrm{opt}, \mu}^f}$
and $\alpha_{t_{{\mathrm{opt}}, p}^f}$, with
\begin{eqnarray*}
t_{\mathrm{opt}, \mu}^f=\inf\bigl\{t_v \dvtx
\mu_{\alpha
_{t_v}}^f\leq b \bigr\} \quad \mbox{and}\quad  t_{\mathrm{opt}, p}^f=
\inf\bigl\{t_v \dvtx  x_{\alpha_{t_v},p}^f\leq b \bigr\},
\end{eqnarray*}
respectively.

Finally, we notice that, usually,
$\mu_{\alpha}^f$ and $x_{\alpha,p}^f$ cannot be derived in a
closed form. Therefore, in order to obtain optimal vaccination
policies, we need to approximate them. The coupling construction
can be used to give a Monte-Carlo based estimation. Suppose, for
simplicity of argument, that $m\leq1$. Fix $n\geq1$, for
$i=1,\ldots,n$, one can simulate a realization $Z^{(i)}$ of $Z$
and $U_j^{(i)}$ of $U_j$, for $j=1,2,\ldots,N^{(i)}(\infty)$,
where $N^{(i)}(\infty)$ is the total number of births in $Z^{(i)}$. For
each $\alpha\in\mathcal{A}^*$, we obtain a realization
$f(Z_{\alpha}^{(i)})$ of $f(Z_{\alpha})$, for $i=1,\ldots,n$. From
these realizations, we estimate $\mu_{\alpha}^f$ and
$x_{\alpha,p}^f$.

\subsection{Time to extinction}
\label{timeextinction} We specialise the preceding results to the case
when evaluation of a
vaccination strategy $\alpha$ is based on the associated distribution of
the time to extinction of the virus in an outbreak. To this end,
for $z\in\mathbb{N}$, we denote by $T_{\alpha,z}$ the time to
extinction of the process $Z_{\alpha}$ when $Z(0)=z$, that is,
\[
T_{\alpha,z}=\inf\bigl\{t\geq0\dvtx  Z_{\alpha}(t)=0\bigr\}.
\]
Thus, $T_{\alpha,z}$ is the maximal time that the infection
survives in the population in an outbreak when the time-dependent
proportion of
immune individuals is given by $\alpha$ and the number of infected
individuals at the beginning of the outbreak is $z$. Now
individuals infect independently of each other, so we have
that
\[
T_{\alpha,z} = \max\bigl\{ T_{\alpha,1}^{(1)},T_{\alpha,1}^{(2)},
\ldots, T_{\alpha,1}^{(z)}\bigr\},
\]
where $T_{\alpha,1}^{(i)}$ are independent
random variables with the same distribution as $T_{\alpha,1}$.
Hence
\[
\mathrm{P}(T_{\alpha,z}\leq t)=\bigl(v_{\alpha}(t)
\bigr)^z,
\]
where
$v_{\alpha}{(t)}=\mathrm{P}(T_{\alpha,1}\leq t)$. Therefore,
to analyze the behaviour of $T_{\alpha,z}$, for any $z$, it is
sufficient to study $T_{\alpha,1}$ through $v_{\alpha}$. From now
on, we denote $T_{\alpha,1}$ by $T_{\alpha}$.

We first use the results of Sections~\ref{sectmoncty} to derive
some continuity and monotonicity properties of the distribution
function $v_{\alpha}$. When every individual is immune, that is, $\alpha
(t)=1$ for all $t > 0$, the infectious
disease does not spread to any susceptible individual and then
the extinction time is given by the survival time of the initial
infected individual. It stands to reason that if there are
non-immune individuals in the population, then it is probable
that the infectious disease takes more time to become extinct. In
the following result, which is an immediate application of Theorem~\ref
{teoa4}(a)
with $f=T$, we show this fact investigating the
behaviour of $v_{\alpha}$ depending on the function
$\alpha$.
%
\begin{corollary}\label{teob1}
Suppose that $\alpha, \alpha' \in\mathcal{A}$ satisfy $\alpha\prec
\alpha'$.
Then $v_{\alpha}(t)\leq
v_{\alpha'}(t)$, for all $t\geq0$.
\end{corollary}

Intuitively, it is clear that the greater the proportion of immune
individuals, the more likely it is that the infectious disease
disappears quickly. Consequently, for any $\alpha\in
\mathcal{A}$, the distribution function $v_{\alpha}$ is bounded
above by $v_{\alpha_1^0}$, the distribution function of the
survival time of the initial infected individual, and bounded
below by $v_{\alpha_0^0}$, which is not necessarily a proper
distribution function. Moreover, we obtain that minor changes in
the proportion of the immune individuals generate minor changes in
the distribution of outbreak duration. The following result is an
immediate application of Theorem~\ref{teoa4}(b), (c) with $f=T$.
%
\begin{corollary}\label{teob2}
\begin{enumerate}[(b)]
\item[{(a)}] Fix
$t>0$. Then, for each $\varepsilon>0$,
\begin{eqnarray*}
\sup_{0\leq u \leq t}\bigl|v_\alpha(u)-v_{\alpha'}(u)\bigr|\leq
\varepsilon,
\end{eqnarray*}
for all $\alpha, \alpha' \in\mathcal{A}$ satisfying
$\|\alpha-\alpha'\|_t \le\delta_{\alpha_0^0,t}(\varepsilon)$.
\item[{(b)}] Fix $t_0\geq0$. Then, for each $\varepsilon>0$,
\begin{eqnarray*}
\sup_{0\leq t < \infty}\bigl|v_\alpha(t)-v_{\alpha'}(t)\bigr|\leq
\varepsilon,
\end{eqnarray*}
for all $\alpha, \alpha' \in\mathcal{A}(c_{\mathrm{inf}},t_0)$
satisfying $\|\alpha-\alpha'\| \le
\delta_{\alpha_{c_{\mathrm{inf}}}^{t_0},\infty}(\varepsilon)$.
\end{enumerate}
\end{corollary}

Finally, we consider the quantiles of $T_{\alpha}$. For $\alpha
\in\mathcal{A}$ and $0 < p < 1$, let
$t_{\alpha,p}=\inf\{t\dvtx v_{\alpha}(t)\ge p\}$
be the quantile of order $p$ of $T_{\alpha}$.
%
\begin{corollary}\label{teob3}
\begin{enumerate}[(b)]
\item[(a)] If $\alpha, \alpha' \in\mathcal{A}$ satisfy $\alpha
\prec\alpha'$, then $t_{\alpha',p} \le t_{\alpha,p}$ for every
$0<p<1$.
\item[(b)] Suppose that $\alpha\in\mathcal{A}$ and
$0<p<1$ are such that $t_{\alpha,p}<\infty$ and $v_\alpha$ is
continuous and strictly increasing at $t_{\alpha,p}$. Then
$\lim_{n \to\infty}t_{\alpha_n,p}=t_{\alpha,p}$, for any sequence
$\{\alpha_n\}$ in $\mathcal{A}$ satisfying $\lim_{n \to\infty}
\alpha_n =\alpha$.
\end{enumerate}%
\end{corollary}
\begin{pf}
\begin{enumerate}[(b)]
\item[(a)] The result follows directly from
Theorem~\ref{teoa5}(a), on setting $f=T$.

\item[(b)] Let $t=t_{\alpha,p}+1$ and $f=\min\{ T,t\}$, so $f \in
\mathcal{P}_t$.
The conditions on $t_{\alpha,p}$ and $v_\alpha$ ensure that $t_{\alpha
,p}=x_{\alpha, p}^f$ for
all $\alpha\in\mathcal{A}$. The result then follows immediately from
Theorem~\ref{teoa5}(b).\qed
\end{enumerate}
\noqed
\end{pf}

Corollary~\ref{teob3} can be extended to a family of vaccination
processes with a continuous index set; cf. Remark~\ref{rr3}(b). In
order to apply Corollary~\ref{teob3}, we
need to determine conditions which guarantee that $v_{\alpha}$ is
both continuous and strictly increasing.
%
\begin{theorem}\label{teob4}
Suppose that the lifetime random variable $I$ is continuous. Then,
for any $\alpha\in\mathcal{A}$, $v_{\alpha}$ is a continuous
distribution function.
\end{theorem}
\begin{pf}
Let $B_0=0$ and, for $n=1,2,\dots\,$, let $B_n$ denote the time of
the $n$th birth in $Z$, with the convention that $B_n=\infty$ if
$N(\infty)<n$. For $n=0,1,\dots,N(\infty)$, let $I_n$ and
$D_n=B_n+I_n$ denote respectively, the lifetime and time of death
of the $n$th individual born in $Z$. Let
$\mathcal{D}=\{D_0,D_1,\dots,D_{N(\infty)}\}$ denote the random
set of all death-times in $Z$. Observe that, for any $t>0$ and any
$\alpha\in\mathcal{A}$, $T_{\alpha}=t$ only
if $t \in\mathcal{D}$. Thus, it is sufficient to show that $\mathrm
{P}(t \in\mathcal{D})=0$ for any $t>0$.

Fix $t>0$ and define $D_n=\infty$ for $n>N(\infty)$. Then, since
$\mathrm{P}(N(t)<\infty)=1$,
%
\begin{equation}
\label{PDt} \mathrm{P}(t \in\mathcal{D})=\mathrm{P} \Biggl(\bigcup
_{n=0}^{\infty}\{ D_n=t\} \Biggr)\le\sum
_{n=0}^{\infty}\mathrm{P}(D_n=t).
\end{equation}
Further, for $n=0,1,\dots\,$,
\begin{eqnarray*}
\label{PDnt} \mathrm{P}(D_n=t)&=&\mathrm{P}\bigl(N(t) \ge n\bigr)
\mathrm{P}\bigl(D_n=t|N(t)\ge n\bigr)
\\
&=&\mathrm{P}\bigl(N(t) \ge n\bigr)\mathrm{E}_{B_n|N(t)\ge n}\bigl[\mathrm {P}
\bigl(D_n=t|B_n,N(t)\ge n\bigr)\bigr]
\\
&=&\mathrm{P}\bigl(N(t) \ge n\bigr)\mathrm{E}_{B_n|N(t)\ge n}\bigl[\mathrm {P}
\bigl(I_n=t-B_n|B_n,N(t)\ge n\bigr)\bigr]
\\
&=&\mathrm{P}\bigl(N(t) \ge n\bigr)\mathrm{E}_{B_n|N(t)\ge n}\bigl[\mathrm
{P}(I_n=t-B_n)\bigr]
\\
&=&0,
\end{eqnarray*}
since $I_n$ is independent of both $B_n$ and $\{N(t)\ge n\}$, and $I$
is continuous.
It then follows from (\ref{PDt}) that $\mathrm{P}(t \in\mathcal
{D})=0$, which completes the proof.
\end{pf}

We notice that under weak conditions, the function $v_\alpha$ is
strictly increasing. Indeed, let $R$ be the number of points of
$\xi$ in $[0,I]$, so $R$ is a random variable giving the number of
offspring of a typical individual in the CMJ branching process
$Z$. Suppose that $\mathrm{P}(R=0)>0$ and that $I|R=0$ is an
absolutely continuous random variable, having density $f_{I|R=0}$
satisfying $f_{I|R=0}(t)>0$ for all $t \in(0,\infty)$. Then it
is easily seen that, for any $\alpha\in\mathcal{A}$,
$v_{\alpha}$ is strictly increasing on $(0,\infty)$, since, for
any open interval $(a,b)$ in $(0,\infty)$, the probability that
the initial individual has no offspring and dies in $(a,b)$ is
strictly positive. It is straightforward to give conditions under
which $v_{\alpha}$ is strictly increasing on $(0,\infty)$ when $I$
has bounded support. For example, suppose that $\mathrm{P}(R=0)$ and
$\mathrm{P}(R=1)$ are both strictly positive, and $I|R=0$ and $B|R=1$
are both absolutely continuous with densities that are strictly
positive on $(0,t_I)$, for some $t_I > 0$. Here, $B$ is the age
that a typical individual has his/her first child. Then, given any
interval $(a,b) \subset(0,\infty)$, there exists $n_0 \in
\mathbb{N}$ such that with strictly positive probability (i) each
of the first $n_0$ individuals in $Z$ has precisely one child,
(ii) the ($n_0+1$)th individual in $Z$ has no children and (iii)
$T \in(a,b)$. It then follows that $\mathrm{P} (T_{\alpha} \in
(a,b) )>0$, provided $\alpha(t)<1$ for all $t>0$.
%
\begin{figure}

\includegraphics{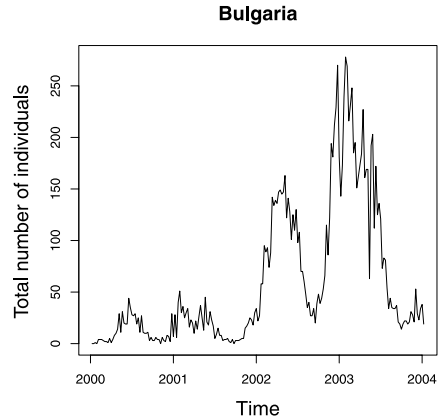}

\caption{Numbers of new infected individuals weekly reported.} \label{f1}
\end{figure}
\section{Illustrative example: Analyzing the control measures for mumps
in Bulgaria}
\label{illus}
As an illustration of how to apply our theoretical results and to
show their usefulness, we analyze a mumps data set
from Bulgaria. In Bulgaria, an increasing number of new cases
of individuals infected with mumps has been observed in recent
years (see Figure~\ref{f1}). This may be a result of a poor
immunization of birth cohorts 1982--1992 (see Kojouharova \textit{et al.} \citep{euro}).
In such a situation, it is necessary to provide supplementary
doses of mumps, measles and rubella (MMR) vaccine targeted at
those cohorts in order to shorten the duration of the outbreaks.
Thus our objective is to determine, using the observed data,
optimal vaccination levels based on the time to extinction that
guarantee, with a high probability, that the outbreak
durations will be less than some suitable bound. {As an example,
we determine} the percentage of the target cohort that must be
vaccinated to guarantee that only primary and first-generation
cases will be observed in at least 90\% of outbreaks.

In order to apply our results, we model the spread of mumps
by a CMJ branching process. This is reasonable since mumps is an infectious
disease which follows the SEIR scheme, and in general, the early
stages of outbreaks following this scheme can be approximated by a
CMJ branching process. Although this is the general situation, a
deeper discussion is needed in the case of mumps. This disease
concerns predominantly young people in schools and universities,
which means small separate populations and population-dependent
propagation. Hence, the approximation of mumps outbreaks in these
populations by CMJ
processes is valid only when outbreaks
are very short, which is the case for the outbreaks we study as we
show later.
%
\begin{figure}

\includegraphics{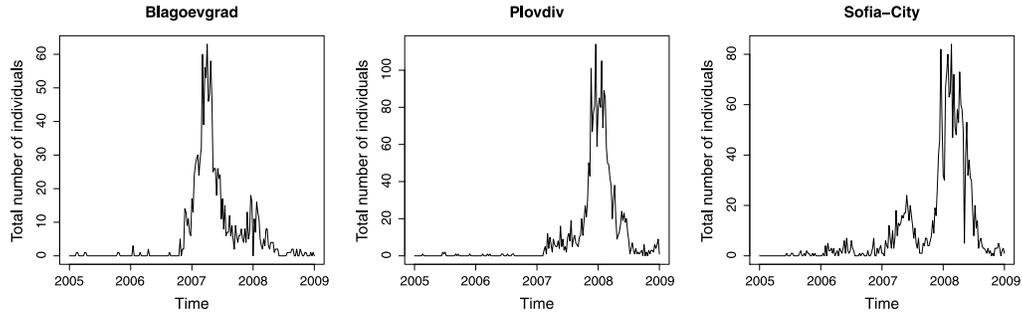}

\caption{Numbers of new infected individuals per week for the
provinces of Bulgaria with the highest incidence of mumps.}
\label{f1_bis}
\end{figure}

The data we analyze (reported by the Bulgarian Ministry of Health)
are the total number of new cases of infected individuals with
mumps observed weekly in each province of Bulgaria from 2005 to
2008, whose birth cohorts were poorly immunized. Notice that we do
not observe outbreak durations, so, first, we describe the
procedure to derive the outbreak durations from these data. Then,
taking into account the main features of mumps transmission, we
select an appropriate general branching process to describe the
evolution of infected individuals in an outbreak and estimate its
main parameters from the data set. Finally, once the model is
fitted, we propose optimal vaccination levels based on the
quantiles of the outbreak duration.

\subsection{Deriving the outbreak duration}

\begin{figure}

\includegraphics{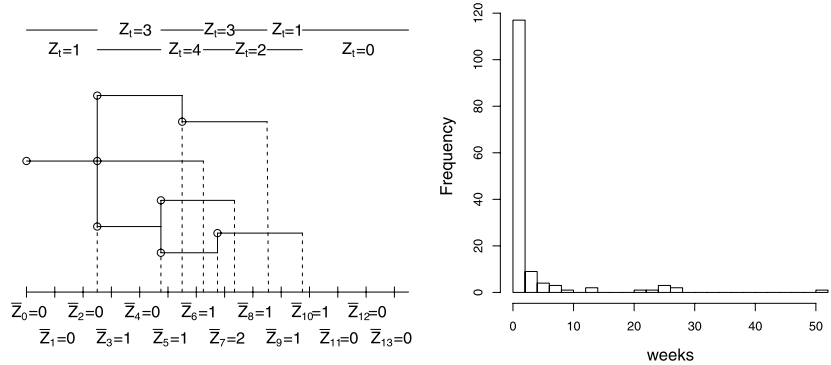}

\caption{Left: Schematic representation of an outbreak. $Z_t$
denotes the underlying branching process and $\bar Z_n$ the number
of new cases in the $n$th week. Right: Durations for outbreaks
started with one infected individual.} \label{zbar}
\end{figure}

Our first task is to determine the behaviour of mumps outbreak
durations in Bulgaria from 2005 to 2008, since our optimal
vaccination level is based on outbreak duration. However, outbreak
durations have not been registered; only the total number of new
cases of infected individuals with mumps in each province has been
observed (see Figure~\ref{f1_bis}). Thus, instead, we derive the
outbreak durations from this data set, taking into account the
main features of mumps transmission. Mumps is a viral infectious
disease of humans and spreads from person to person through the
air. The period between {someone being transmitted mumps and that
person first showing symptoms of mumps is called the incubation
period for mumps}. This incubation period can be 12 to 25 days and
the average is 16 to 18 days. The infectious period {(i.e., when an
individual is able to transmit the mumps virus to others)} starts
about 2 days before the onset of symptoms and usually, an
individual with mumps symptoms is immediately isolated from the
population (see \url{http://kidshealth.org}). In view of the range
of the incubation period, we consider that an outbreak is formed
by the cases that appear in a province in a sequence of weeks with
no more than three consecutive weeks without cases. That is, when
we observe more than three weeks without cases we consider that
the outbreak has become extinct, with the next outbreak starting
in the first subsequent week in which there is at least one new
case. Applying this procedure for each province, we have obtained
262 outbreaks. The left plot in Figure~\ref{zbar} could represent
one such outbreak initiated by one infected individual. In this
schematic representation, we have considered that the infectious
period is \emph{negligible} due to the fact that infected
individuals are immediately isolated when they show symptoms. The
variable $Z_t$ denotes the underlying branching process, {which is
not} observed. The segments over/under $Z_t$ indicates the lengths
of time for which $Z_t$ takes the corresponding values. The tick
marks on the axis represent weeks, and $\bar Z_n$ the number of
new cases observed during the $n$th week. Indeed, $\bar Z_n$, $n
\geq0$, are the variables that {are} observed. In this context,
by outbreak duration we mean the time elapsing between the
appearance of the first case until isolation of the last one, that
is the time to extinction of the branching process minus the
incubation period of the first individual. Thus, a more accurate
way to approximate outbreak duration from the observed data is by
the total number of weeks until extinction of the virus (giving an
error, due to discretization, of at most one week), yielding seven
weeks in the outbreak of Figure~3 (left).
%

For each of the 262 outbreaks, we calculated the total number
of weeks until extinction of the virus (and, also, the outbreak size,
i.e., total
number of infected individuals). We noticed that
the behaviour of these outbreak durations depends on the initial
number of infected individuals. Hence, we have considered only
those outbreaks which started with one infected individual, a
total of 144. We checked that both outbreak duration and outbreak size
were homogeneous between provinces (Kruskal--Wallis test: $p$-values
0.4763 and
0.4782, resp.) and consequently assumed that disease
propagation in the different provinces are independent
replications of the same process. Thus, the right plot in
Figure~\ref{zbar} shows the histogram of outbreak durations for
all 144 outbreaks started with one infected individual. We
observe two different groups, outbreaks for which their duration
is less than 10 weeks (comprising 134 outbreaks) and another group
where the outbreak duration is greater than 10 weeks (comprising
the remaining ten outbreaks). Possibly, this happens because some
cases observed in a week could not come from cases of previous
weeks, and then new outbreaks could have appeared overlapping
in time. Hence, we consider that the outbreaks corresponding to
durations of this last group may have been initiated no more than
10 weeks before. Thus, outbreak durations greater than 10 weeks
have been removed from our study, and only durations less than 10
weeks have been considered in order not to overestimate the
duration of the outbreaks. Nevertheless, an outbreak with apparent
duration less than 10 weeks could actually be the superposition of
two or more separate outbreaks, but we cannot determine this.
%
\begin{figure}

\includegraphics{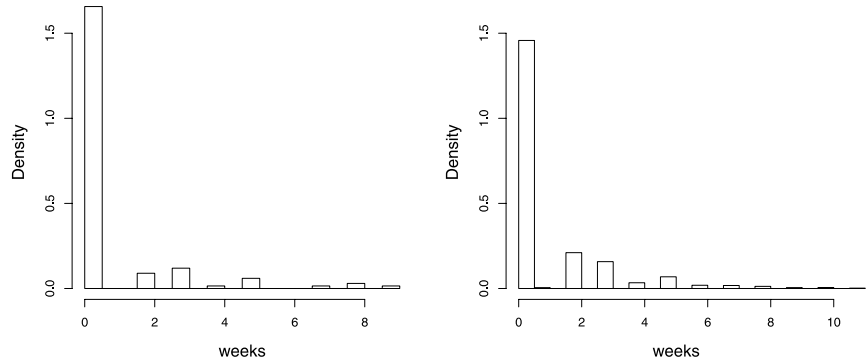}

\caption{Left: Durations for outbreaks started with one infected
individual without overlapping. Right: Simulated durations from a
BHBP for outbreaks started with one infected individual.}
\label{f2}
\end{figure}

The left plot of Figure~\ref{f2} shows the durations of the 134
outbreaks considered. We notice that 83\% of these outbreaks have
only one infected individual, so their outbreak durations are 0. {
The remaining 17\% of outbreaks seem to have a cyclical behaviour
with period given by the mean of the incubation period
(approximately 2.5 weeks)}.

\subsection{Modelling mumps transmission}
\label{mumps}

As noted above, mumps is a contagious disease of humans that is
spread from person to person through the air. The most common
method of transmission is through coughing or sneezing, which can
spread droplets of saliva and mucus infected with the mumps virus.
Hence, when an infected person coughs or sneezes, the droplets
atomize and can enter the eyes, nose, or mouth of another person.
Following mumps transmission, a person does not immediately become
sick. Once the virus enters the body, it travels to the back of
the throat, nose and lymph glands in the neck, where it begins to
multiply. As indicated previously, this period between mumps
transmission and the beginning of mumps symptoms is the incubation
period for mumps. People who have mumps are most contagious from 2
days before symptoms begin to 6 days after they end and
transmission may occur at anytime in that period. Since an
individual with mumps symptoms is immediately isolated from the
population, the infectious period is very short in comparison with
the incubation period, so, as indicated previously, we assume that
transmission occurs only
at the end point of an individual's incubation period. This assumption
simplifies the
mathematical model and does not influence strongly outbreak
duration.
As the end of the incubation period means that an individual's viral load
has reached a given threshold to produce clinical signs,
we assume that the mean number of individuals infected by an
infected individual is constant and does not depend on the
length of his/her incubation period.

An earlier analysis of these mumps data using
Bienaym{\'e}--Galton--Watson branching processes is given in Angelov
and Slavtchova--Bojkova \citep{angelov2012}. However, the above
observations imply that the Bellman--Harris branching process
(BHBP) (see Athreya and Ney \citep{athreya}) is a more appropriate
model for mumps transmission and indeed it provides an improved
fit to these data. Recall that a BHBP is a CMJ branching process,
in which an individual reproduces only at the end of his/her
life-time, according to an offspring law which is the same for all
the individuals. In the epidemiological context, age is the
incubation period and the reproduction law is the contagion
distribution.

Next, we describe the incubation period and contagion
distributions used to model mumps transmission in each outbreak in
Bulgaria by means of the same BHBP (recall that we did not find
any difference in the behaviour of the outbreaks in different
provinces). We assume that the incubation period $I$ follows a
gamma distribution, with shape parameter $r>0$ and rate
$\gamma>0$, so $I$ has mean $r\gamma^{-1}$ and probability density
function
\begin{eqnarray*}
f_I(u)=\frac{\gamma^ru^{r-1}\exp(-\gamma u)}{\Gamma(r)},\qquad u>0,
\end{eqnarray*}
where $\Gamma$ is the gamma function, and that the contagion
distribution follows a Poisson distribution with mean $m$. These
distributions are appropriate for the incubation period and the
number of infections, respectively (see, e.g., Daley and
Gani \citep{daleygani}, Farrington and Grant \citep{fg},
Farrington \textit{et al.} \citep{fkg} or Mode and
Sleeman \citep{modesleemam}). Intuitively, $m$, the mean number of
individuals infected by an infected individual, represents the
power of the virus. Taking into account that the incubation period
is estimated between 12 and 25 days and the average is 16 to 18
days, we consider the gamma distribution with mean 17 and $r=50$,
{which implies} that the incubation period in $98.7\%$ of
individuals is between 12 and 25 days. To estimate $m$, we consider
the maximum likelihood estimator (MLE) based on the total number
of births in independent extinct realisations of a BHBP. The total
number of births in a BHBP has the same distribution as that in a
Bienaym\'{e}--Galton--Watson branching process with the same offspring
distribution. In our application, the offspring distribution is
Poisson and it follows that the total number of births $N(\infty)$
(excluding the initial $a$ individuals) follows a Borel--Tanner
distribution with probability mass function
\[
\mathrm{P}\bigl(N(\infty)=k\bigr)=\frac{a m^k (a+k)^{k-1} \mathrm
{e}^{-(a+k)m}}{k!}, \qquad k=0,1,\dots.
\]
(Note that, for $l=1,2,\dots\,$, the mean number of births in the
$l$th generation is $a m^l$, so the expectation of this
 Borel--Tanner distribution is
$E[N(\infty)]=a(m+m^2+\dots)=am(1-m)^{-1}$, when~$m<1$.) It follows that the
MLE of the offspring mean $m$, based on $L$ independent
realisations, is given by $\hat m = (\sum_{i=1}^L
n^{(i)})(\sum_{i=1}^L a^{(i)} + n^{(i)})^{-1}$, where, for
$i=1,2,\ldots,L$, $a^{(i)}$ and $n^{(i)}$ are respectively the
initial number of individuals and the total number of births in
the $i$th realisation (for details see Farrington \textit{et
al.} \cite{fkg}). In our case $L=134$, $\sum_{i=1}^L a^{(i)}=134$
and $\sum_{i=1}^L n^{(i)}=62$, whence $\hat m = 0.3163$. Note that
inference based on duration of outbreaks is less sensitive to
underreporting than that based on the total number of births.
However, estimating the offspring law based on the time to
extinction of each outbreak turns into a difficult problem in
branching processes {theory, even for the simplest model} (see, e.g., Farrington \textit{et al.} \citep{fkg}).

Applying the general theory of branching processes, since the
estimated value of $m$ is less than~1, we deduce that mumps
transmission can still occur in Bulgaria, but such spread cannot
lead to a large-scale epidemic. This fact is consistent with the
Figures~\ref{f1} and \ref{f1_bis}. Although the epidemic becomes
extinct, it can have different levels of severity. One measure of
severity is the mean size of {an} outbreak, excluding the initial
case, viz. $m(1-m)^{-1}$, which in our case is estimated by 0.463.
However, we are concerned with the problem of how to shorten
outbreak durations by vaccination. To this end, we analyze the
random variable $T_{\alpha_{c_{\mathrm{inf}}}^0}$, the time to extinction
of a BHBP with incubation period and contagion distributions as
described above. {Note that $c_{\mathrm{inf}}=0$, as $m\le1$, so here
$T_{\alpha_{c_{\mathrm{inf}}}^0}$ is the extinction time when there is no
supplementary vaccination.} {The }variable
$T_{\alpha_{c_{\mathrm{inf}}}^0}$ includes the incubation period of the
initial individual, which is not observed in practice. Thus, from
now on, we use the random variable
$\widetilde{T}_{\alpha_{c_{\mathrm{inf}}}^0}$, the difference between
$T_{\alpha_{c_{\mathrm{inf}}}^0}$ and the incubation period of the initial
individual (i.e., the definition of outbreak duration given in the
previous subsection) to model mumps outbreak duration in Bulgaria.
The right plot in Figure~\ref{f2} shows a histogram of $10\,000$
simulated durations of outbreaks (rounded up to the nearest
integer), each initiated by one infected individual and modelled
by a BHBP with the above parameters. We notice that in 72.9\% of
these simulated outbreaks the initial infected individual does not
infect any new individual {(recall 83\% for real data)}. Moreover,
{the simulated outbreak durations show the same cyclical behaviour
as seen in the real data.}

Comparing real and simulated durations, we deduce that mumps
outbreak durations in Bulgaria can be modelled by the variable
$\widetilde{T}_{\alpha_{c_{\mathrm{inf}}}^0}$ (Pearson's chi-squared test:
$p$-value 0.2951, grouping the tail for values greater than 8).

\subsection{Determining the optimal vaccination levels}

Once we have fitted the model, in order to apply our theoretical
results we have assumed that the proportion of immune individuals
is constant with time, since, generally, vaccination is applied
when an individual is a child and the disease spreads when he/she
is a teenager. In the particular case of supplementary vaccination
for Bulgarian mumps, for simplicity we assume that this
vaccination process occurs simultaneously across the country (e.g., in secondary schools at the same specific time). To
determine the optimal vaccination levels, we denote by
$\widetilde{T}_{\alpha_{c}^0}$ the difference between
$T_{\alpha_{c}^0}$ and the incubation period of the initial
individual, when the proportion of immune individuals in the
population is $c$, with $0\leq c\leq1$. In the same way as was
proved for $T_{\alpha_{c}^0}$ (see Corollary~\ref{teob3}), we
deduce that $\widetilde{T}_{\alpha_{c}^0}$ has the same quantile
properties depending on $c$ as $T_{\alpha_{c}^0}$ (notice that
$\widetilde{T}_{\alpha_{c}^0}$ is monotonically decreasing with
pruning). Therefore, next we propose vaccination policies based on
the quantiles of $\widetilde{T}_{\alpha_{c}^0}$, with $0\leq c\leq
1$. Specifically, for fixed $p$ and $t$, with $0 < p < 1$ and $t
> 0$, we seek vaccination policies which guarantee that the
mumps virus becomes extinct in each outbreak, with probability
greater than or equal to $p$, not later than time $t$ after the
outbreak has been detected with $z$ initial infected individuals,
that is\looseness=-1
\begin{eqnarray*}
c_{\mathrm{opt}}=c_{\mathrm{opt}}(z,p,t)= \inf\bigl\{c\dvtx  0\leq c\leq
1,{x}^{\widetilde{T}}_{\alpha_c^0,p^{1/z}}\leq t\bigr\},
\end{eqnarray*}\looseness=0
where ${x}^{\widetilde{T}}_{\alpha_c^0,p^{1/z}}$ denotes the
quantile of order $p^{1/z}$ of the variable
$\widetilde{T}_{\alpha_c^0}$.\vadjust{\goodbreak}
%
\begin{figure}

\includegraphics{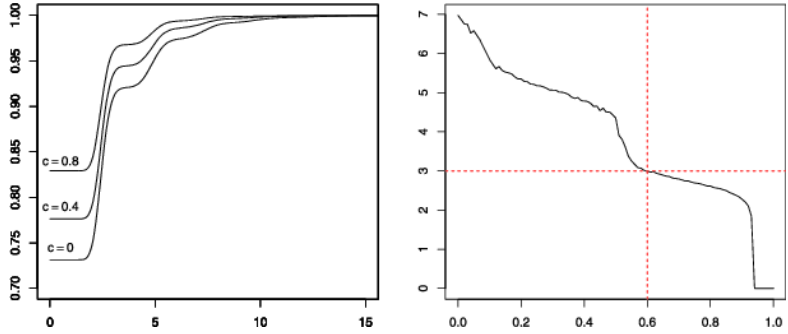}

\caption{Left: Behaviour of the distribution function of
$\widetilde{T}_{\alpha_{c}^0}$ for $c=0,0.4,0.8$. Right: Behaviour
of ${x}^{\widetilde{T}}_{\alpha_c^0,0.9^{1/5}}$ depending on $c$,
with $0\leq c\leq1$.} \label{f3}
\end{figure}

As an illustration, we take $z=5$, $p=0.9$ and $t=3$, being the
time measured in weeks. First we justify these values. Consider
the value of $z$. Since the number of infected individuals at the
beginning of an outbreak is unknown, we bound it by the greatest
number of individuals infected by one infected individual. Taking
into account that the contagion distribution is Poisson and the
estimate of $m$, we obtain the upper bound to be $5$, and
therefore we take $z=5$. Moreover, we select $t=3$, which, taking
into account the features of the incubation period, guarantees
that only primary and first-generation cases will be observed.
Since in our situation the estimated value of $m$ is less than~1,
to approximate $c_{\mathrm{opt}}$, we need to obtain the empirical
distribution of $\widetilde{T}_{\alpha_{c}^0}$, for $0\leq c\leq
1$, using the Monte-Carlo method described in
Section~\ref{optimal}. To this end, for each $c=0.01k$, with
$k=0,\ldots,100$, $100\,000$ processes have been simulated and
their duration calculated. The left plot in Figure~\ref{f3} shows
the behaviour of the empirical distribution function of
$\widetilde{T}_{\alpha_{c}^0}$ for several values of $c$. Notice
that as $c$ increases, the outbreak duration decreases in a
continuous way, in accordance with Corollaries \ref{teob1}
and \ref{teob2}. The right plot in Figure~\ref{f3} shows the
behaviour of ${x}^{\widetilde{T}}_{\alpha_c^0,0.9^{1/5}}$
depending on $c$, which is in accordance with
Corollary~\ref{teob3}. Since
${x}^{\widetilde{T}}_{\alpha_{c_{\mathrm{inf}}}^0,0.9^{1/5}}=6.97$, our
model estimates that the duration of 90\% of outbreaks in Bulgaria
is less than 6.97 weeks, if vaccination is not applied (in our
real data 97\% of outbreaks have durations less than 6 weeks). In
order to shorten the outbreak duration, from our study, we deduce
that $c_{\mathrm{opt}}(5,0.9,3)=0.6$ (see right plot in Figure~\ref{f3}).
Therefore, vaccinating a proportion of 60\% of susceptible
individuals in the target cohort, guarantees that in at least 90\%
of outbreaks of mumps in Bulgaria only primary and
first-generation cases will be observed after the vaccination.
Finally, we notice that $c_{\mathrm{opt}}(5,0.9,0)=0.94$, that is, to
guarantee that at least the
90\% of outbreaks do not spread after vaccination, the vaccination
level should
be 94\% of susceptible individuals in the target cohort.
%
\begin{table}
\tablewidth=\textwidth
\tabcolsep=0pt
\caption{Sensitivity analysis on the mean and shape parameter of
the gamma incubation distribution} \label{sa}
\begin{tabular*}{\textwidth}{@{\extracolsep{\fill}}d{2.1}llllll@{}}
\hline
\multicolumn{1}{l}{Mean}& &\multicolumn{5}{c}{Shape
parameter $r$}\\
[-5pt]
& &\multicolumn{5}{c}{\hrulefill}\\
 &  &\multicolumn{1}{l}{30} & \multicolumn{1}{l}{40}&
\multicolumn{1}{l}{50} & \multicolumn{1}{l}{60}& \multicolumn{1}{l}{70} \\ \hline
16 &\%
Coverage & 92.2&95.3&97.1&98.8&98.8\\
& $c_{\mathrm{opt}}(5,0.9,3)$&0.60&0.57&0.56&0.54&0.54\\
16.5 &\%
Coverage & 93.0&96.6&98.1&98.9&99.4\\
& $c_{\mathrm{opt}}(5,0.9,3)$&0.63&0.60&0.58&0.56&0.55\\
17 &\%
Coverage & 94.9&95.5&98.7&99.3&99.6\\
& $c_{\mathrm{opt}}(5,0.9,3)$&0.66&0.64&0.60&0.58&0.57\\
17.5 &\%
Coverage & 95.4&97.9&99.0&99.5&99.8\\
& $c_{\mathrm{opt}}(5,0.9,3)$&0.70&0.67&0.65&0.62&0.61\\
18 &\%
Coverage & 95.3&97.8&99.0&99.5&99.8\\
& $c_{\mathrm{opt}}(5,0.9,3)$&0.73&0.71&0.68&0.65&0.64\\
\hline
\end{tabular*}
\end{table}

The parameters of the gamma distribution used to model the
incubation period have been derived from knowledge of mumps
transmission rather than estimated from data. Thus we have
performed a sensitivity analysis of their influence on the optimal
vaccination level. {We have considered gamma distributions with
mean and shape parameter $r$ taking values in a grid (giving
different probabilities for the incubation period belonging to
range 12--25, which we denote as percentages of coverage),
yielding the results shown in Table~\ref{sa}. One can observe that
increasing the mean (holding $r$ fixed) clearly increases the
duration of the epidemic leading to higher values of $c_{\mathrm{opt}}$.
Moreover, increasing the shape parameter $r$ (holding the mean
fixed) decreases the variance of lifetimes and hence also the
chance of long outbreak duration, leading to lower values of
$c_{\mathrm{opt}}$. The optimal vaccination level $c_{\mathrm{opt}}(5,0.9,3)$ is
fairly stable in the vicinity of the chosen values of 17 and 50
for the mean and shape parameter $r$, respectively.}
%
\begin{remark}
From a computational point of view it is interesting to note that
to find optimal vaccination policies, the simulation method based
on \emph{pruning}, described at the end of Section~\ref{optimal},
has proved to be at least 17\% faster than those in Gonz\'alez \textit{et
al.} \citep{gmb1,gmb2}, {which are also simulation-based
methods but work directly with the distribution of the extinction
time. For the BHBP there exist other methods to approximate the
distribution function of the time to extinction based on solving
numerically an associated integral equation (see Mart{\'\i}nez and
Slavtchova-Bojkova \citep{rm}, which includes comparison with
simulation-based methods). Unlike the latter approach, the
{Monte-Carlo} method proposed in Section~\ref{optimal} is easily
extended to time-dependent vaccination processes.} All the
computations and simulations have been made with the statistical
computing and graphics language and environment $\mathbf{R}$
(``GNU S'', see \cite{r}).
\end{remark}

\section{Concluding comments}
\label{conc} The coupled pruning technique for proving
monotonicity and continuity properties of functions defined on CMJ
branching processes depending on the vaccination function $\alpha$
is both simple and powerful. It is clear that the proofs
generalise easily to more general branching processes, such as
multitype CMJ branching processes, time-inhomogeneous branching
processes and branching processes in a random environment. The
function $\alpha$ does not have to represent vaccination. It
could represent any control of disease propagation that has the
effect of reducing either the number of susceptibles or the
probability that a contacted susceptible becomes infected.
However, for the coupled pruning technique to work it is necessary
that, in the branching process setting, the control affects only
the probability that a birth is aborted and not the intrinsic
reproduction law of the branching process. Thus, for example, the
method cannot be applied to density-dependent processes, such as
population size dependent branching processes, if the density
dependence relates to the size of the unvaccinated population
rather than the total population size.

Given that the results in the Bulgarian mumps illustration are
based on simulation alone, it may seem more appropriate to use an
epidemic model rather than a branching process that approximates
such a model. However, there are several advantages in using the
simpler branching process formulation. First, branching process
models can be fitted directly to the data more easily; in
particular they do not require knowledge of the size of the
population in which the outbreaks are occurring. Second, the
coupled pruning technique enables the monotonicity and continuity
properties pertaining to vaccination functions to be proved
easily. Third, the coupled pruning technique yields an
associated {Monte-Carlo} method for determining optimal
vaccination processes. 

The framework for optimal vaccination policies studied in
Section~\ref{optimal} can be extended to include alternative
formulations of optimal policies. For example, one may define a
cost $c(\alpha)$ associated with each vaccination process $\alpha
\in\mathcal{A}$ and then seek vaccination processes from a subset
$\mathcal{A}^*$ of $\mathcal{A}$ which either (i) minimise
$c(\alpha)$ subject to $\mu_\alpha^f \le b$ or (ii) minimise
$\mu_\alpha^f$ subject to $c(\alpha) \le c_0$, where $c_0$ is
specified. Provided the cost function, $c(\alpha)$ is suitably
monotonic and continuous in $\alpha$ and $\mathcal{A}^*$ is
totally ordered, Theorems \ref{teoa1} and \ref{teoa2} imply the
existence of unique such optimal vaccination processes and it
should be possible to extend the {Monte-Carlo} algorithm at the
end of Section~\ref{optimal} to estimate the optimal vaccination
processes. {Optimal vaccination} policies that permit vaccination
costs to be taken into account are especially relevant in animal
vaccination.

\section*{Acknowledgements}
We thank the referees for their careful reading of our paper and
for their constructive comments which have improved its
presentation. The research was partially supported by the
Ministerio de Econom{\'\i}a y Competitividad and the FEDER through
the Plan Nacional de Investigaci\'on Cient{\'\i}fica, Desarrollo e
Innovaci\'on Tecnol\'ogica, grant MTM2012-31235 and
by the appropriated state funds for research allocated to Sofia
University, grant 125/2012, Bulgaria.




\printhistory

\end{document}